\documentclass[10pt]{article}
\usepackage{amssymb}
\usepackage{mathrsfs}
\usepackage{amsfonts}
\usepackage{tikz}

\usepackage{amsmath,amsthm,amsfonts,amssymb,amscd}

\def\ex{\mathrm{ex}}

\allowdisplaybreaks[4]
\newtheorem {Problem} {Problem}[section]
\newtheorem {Theorem} [Problem]{Theorem}
\newtheorem {Lemma}[Problem]{Lemma}

\newtheorem {Corollary}[Problem]{Corollary}
\newenvironment {Proof}{\noindent {\bf Proof.}}{\hfill\ensuremath{\square}}
\newcommand*{\QEDB}{\hfill\ensuremath{\square}}

\begin{document}

\title{On the $A_\alpha$-spectral radius of  graphs without linear forests \thanks{This work is supported by  the National Natural Science Foundation of China (Nos. 12101166, 12101165, 11971311, 12161141003),
and the Hainan Provincial Natural Science Foundation of China (Nos. 120RC453, 120MS002).
\newline \indent Email: mzchen@hainanu.edu.cn, amliu@hainanu.edu.cn, xiaodong@sjtu.edu.cn. \newline  \indent $^{\dagger}$Corresponding author:
Xiao-Dong Zhang (Email: xiaodong@sjtu.edu.cn),}}

\author{ Ming-Zhu Chen, A-Ming Liu,\\
School of Science, Hainan University, Haikou 570228, P. R. China, \\
\and  Xiao-Dong Zhang$^{\dagger}$
\\School of Mathematical Sciences, MOE-LSC, SHL-MAC\\
Shanghai Jiao Tong University,
Shanghai 200240, P. R. China}

\date{}
\maketitle

\begin{abstract}
Let $A(G)$ and $D(G)$ be the adjacency  and  degree matrices of a simple graph $G$ on $n$ vertices, respectively.
The \emph{$A_\alpha$-spectral radius} of $G$ is the largest eigenvalue of $A_\alpha (G)=\alpha D(G)+(1-\alpha)A(G)$ for a  real number $\alpha \in[0,1]$.
In this paper, for $\alpha \in (0,1)$, we obtain a sharp upper bound for the  $A_\alpha$-spectral radius of  graphs on $n$ vertices  without a subgraph isomorphic to a  liner forest  for $n$ large enough and characterize all graphs which attain the upper bound.  As a result,
we completely obtain the maximum signless Laplacian spectral radius of graphs on $n$ vertices without a subgraph isomorphic to a  liner forest  for $n$ large enough. 
\end{abstract}

\section{Introduction}
Let $G$ be an undirected simple graph with vertex set
$V(G)=\{v_1,\dots,v_n\}$ and edge set $E(G)$.
The \emph{adjacency matrix}
$A(G)$ of $G$  is the $n\times n$ matrix $(a_{ij})$, where
$a_{ij}=1$ if $v_i$ is adjacent to $v_j$, and $0$ otherwise.   The  \emph {spectral radius} of $G$, denoted by $\rho(G)$,   is the largest eigenvalue of $A(G)$. The  \emph{signless Laplacian spectral radius} of $G$,  denoted by $q(G)$, is the largest eigenvalue of $Q(G)$, where $Q(G)=A(G)+D(G)$ and $D(G)$ is the degree diagonal matrix.
For $v\in V(G)$,  
the \emph{degree} $d_G(v)$ of $v$  is the number of vertices adjacent to $v$ in $G$.
We write  $d(v)$ for $d_G(v)$  if there is no ambiguity.
For a set of vertices $X$, we use 
$N_G[X]$ to denote $\bigcap_{v\in X}N_G(v)$. 
 For two sets $X,Y\subseteq V(G)$,  we use $e(X)$ to denote the number of edges in $G$ with two ends in $X$ and $e(X,Y)$ to  denote the number of edges in $G$ with one end in $X$ and the other in $Y$.  For two vertex disjoint graphs $G$ and $H$,  we denote by  $G\cup H$ and  $G\nabla H$  the \emph{union} of $G$ and $H$,
and the \emph{join} of $G$ and $H$, respectively.
Denote by $kG$  and $\overline{G}$ the union of $k$ disjoint copies of $G$ and the complement graph of $G$, respectively.  We say that a graph $G$ is \emph{$H$-free} if  it does
not contain a subgraph isomorphic to $H$, i.e., $G$ contains no copy of $H$.
The {\it Tur\'{a}n number}, denoted by $\ex(n,H)$,  of a graph $H $ is the maximum number of edges in  an $H$-free graph on $n$ vertices.
 For a bipartite graph $H$, the bipartite Tur\'{a}n number, denoted by $\ex(m,n; H)$, of a bipartite graph $H$ is the maximum number of edges in an $H$-free bipartite graph $G=(X,Y;E)$ with $|X|=m$, $|Y|=n$ and $m\le n$.

To track the gradual change of $A (G)$ into $Q (G)$,
Nikiforov \cite{Nikiforov2017} proposed  and studied the convex linear combinations $A_\alpha(G)$ of $A (G)$ and
$D (G) $ defined by
$$A_\alpha (G)=\alpha D(G)+(1-\alpha)A(G)$$
for any given real number $\alpha\in [0,1]$.
  Note that $A_0(G)=A(G)$, $2A_{1/2}(G)=Q(G)$, and $A_1(G)=D(G)$. The  \emph{$A_\alpha$-spectral radius} (or \emph{$\alpha$-index}) of $G$ is the largest eigenvalue of $A_\alpha(G)$, denoted by
 $\rho_\alpha(G)$. Clearly, $\rho_0(G)=\rho(G)$ and $2\rho_{1/2}(G)=q(G)$. The eigenvalues of  $A_\alpha(G)$ are intensively investigated in the literature, see \cite{Chen2019,Guo2022,He2021,Lin2019,Nikiforov2017,Tian2021,Yu2019}.

 Note that $A_\alpha(G)$ is nonnegative. 
For $\alpha \in [0,1)$, if $ G$ is connected then by the Perron-Frobenius theorem, we know that $A_\alpha(G)$ has
an eigenvector corresponding to $\rho_\alpha(G)$ with all entries being positive, denoted by $\mathbf x$, also see \cite{Nikiforov2017}.
For a vertex $v \in V (G)$, we will write $x_v$ for the eigenvector entry of $\mathbf x$ corresponding to
$v$. We may normalize $\mathbf x$ so that it has maximum entry equal to $1$.
 If there are multiple such vertices, we choose and fix one of them arbitrarily
among them. In addition, if $G$ is connected and $H$ is a proper subgraph of $G$, then $\rho_\alpha(G)>\rho_\alpha(H)$.
Further, by eigenequations of $A_\alpha(G)$ on any vertex $u\in V(G)$,
$$\rho_\alpha(G)x_u=\alpha d(u) x_u+(1-\alpha)\sum_{v\in N(u)} x_v;$$
and by Rayleigh's principle (for example, see \cite{Nikiforov2017}),
$$\rho_\alpha (G)=\max_{||{\bf x}||_2=1}\sum_{uv\in E(G)}(\alpha x_u^2+2(1-\alpha)x_u+\alpha x_v^2).$$

In this paper, for $\alpha \in (0,1)$, we determine the maximum $A_\alpha$-spectral radius
of graphs on $n$ vertices which do not contain a subgraph isomorphic to a linear forest for  $n$ large enough by an unified approach.
Denote by $S_{n,p}=K_p\nabla \overline{K}_{n-p}$, $S^+_{n,p}=K_p\nabla (\overline{K}_{n-p-2}\cup K_2)$,
 and $F_{n,p}=K_{p}\nabla ( t K_2\cup rK_1)$, where $n-p=2t+r$ and $0\leq r<2$. The main result of this paper can be stated as follows.

 \begin{Theorem}\label{linear forest}\label{thm}
 Let  $F_{\ell}=\cup_{i=1}^{\ell} P_{k_i}$ be a linear forest with   $\ell\geq1$ and
$k_1\geq \cdots \geq k_{\ell}\geq2$ and $ p=\sum _{i=1}^{\ell} \lfloor\frac{k_i}{2}\rfloor-1\geq1$.  Suppose that $G$ is  an $F_{\ell}$-free  graph  of sufficiently large order $n$. If $\alpha \in (0,1)$, then the following holds.\\
(i) If there  exists  at least one even $k_i$, then $$\rho_\alpha(G)\leq  \frac{\alpha n+p-1+\sqrt{(\alpha n+p-1)^2-4p(2\alpha-1)n-4p(p-(p+1)\alpha)}}{2},$$ and the equality holds if and only if $G=S_{n,p}$;\\
(ii) If  all $k_i$  are  odd and there exists one $k_i>3$,  then $\rho_\alpha(G)$ is no more than the largest root of $f(x)=0$,
and  the equality holds if and only if $G=S^+_{n,p}$,
where $f(x)=x^3 - (p + n\alpha + p\alpha)x^2 + (p - np + n\alpha - 2 p\alpha + p^2 + np\alpha^2 + 2np\alpha - 1)x + p^2\alpha^2 - 2p + p^3\alpha^2 + np + 5 p\alpha - 2 \alpha^2 + p^2\alpha - p^3\alpha - p^2 + np^2\alpha - 2np^2\alpha^2 - 2np\alpha$.\\
(iii) If $k_1=\cdots=k_{\ell}= 3$, then $\rho_\alpha(G)\leq t_\alpha(n,p)$, and  the equality holds if and only if   $G=F_{n,p}$,
where
 $t_\alpha(n,p)=\frac{p+\alpha n+\sqrt{(\alpha n-3p)^2+4(p-\alpha)n+4(p+2)(\alpha p-3p+5)}}{2}$  if $n-p$ is even and
$t_\alpha(n,p)$ is the largest root of $f(x)-p(\alpha-1)^2(n-p-3)=0$  if $n-p$ is odd.
\end{Theorem}

 Let $\alpha=\frac{1}{2}$. Then we can easily obtain the following corollary  which gives a complete solution to the maximum signless Laplacian spectral radius of graphs without a subgraph isomorphic to a linear forest, which improves the result in \cite{Chen2020}.
 \begin{Corollary}\label{Q spectral radius for all linear forest}
Let  $F_{\ell}=\bigcup_{i=1}^{\ell} P_{k_i}$ be a linear forest with   $\ell\geq 1$ and
$k_1\geq \cdots \geq k_{\ell}\geq2$ and $ p=\sum _{i=1}^{\ell} \lfloor\frac{k_i}{2}\rfloor-1\geq 1$.  Suppose that $G$ is  an $F_{\ell}$-free  graph  of sufficiently large order $n$. Then the following holds.\\
(i) If there  exists  at least one even $k_i$, then $q(G)\leq \frac{2n+2p-2+\sqrt{(n+2p-2)^2-8p(p-1)}}{2}$, and the  equality holds if and only if $G=S_{n,p}$;\\
(ii) If  all $k_i$  are  odd and there exists one $k_i>3$,  then $q(G)$ is no more than the largest root of $g(x)=0$, where $g(x)=x^3 - (n + 3p)x^2 + (4p^2 + np + 2n - 4)x - 2p^2(p+1)  =0$,
and  the equality holds if and only if $G=S^+_{n,p}$;\\
(iii) If $k_1=\cdots=k_{\ell}= 3$, then $q(G)\leq 2t_{1/2}(n,p)$, and  the equality holds if and only if   $G=F_{n,p}$,
where
 $2t_{1/2}(n,p)=\frac{n+2p+\sqrt{(n+2p-4)^2-8p(p-2)}}{2}$  if $n-p$ is even and
$2t_{1/2}(n,p)$ is the largest root of $g(x)-2p(n - p - 3)=0$  if $n-p$ is odd.
\end{Corollary}


The rest of this paper is organized as follows. In Section~2, we present some known and necessary results.
In Section~3, we give the proof of Theorem~\ref{thm}.

\section{Preliminary}
In 1959, Erd\H{o}s and Gallai in \cite{Erdos1959} proved the following key result, which opens a new subject for extremal graph theory.

\begin{Lemma}\label{lem2.1}\cite{Erdos1959}(Erd\H{o}s-Gallai theorem)
Let $G$ be a   graph of order $n$.   If  $G$ is $P_k$-free, then $e(G)\leq(k-2)n/2$ with equality if and only if $G$ is a  union of disjoint copies of $K_{k-1}$'s.
\end{Lemma}
%
In 2013, Lidick\'{y} \cite{Lidicky2013} extended Erd\H{o}s-Gallai theorem to a linear forest with at least a path of order not equal to $3$.
\begin{Lemma}\label{lem2.2}\cite{Lidicky2013}
 Let $F_{\ell}=\bigcup_{i=1}^{\ell} P_{k_i}$ be a linear forest with   ${\ell}\geq 2$ and
$k_1\geq \cdots \geq k_{\ell}\geq2$ and $ p=\sum _{i=1}^{\ell} \lfloor k_i/2 \rfloor-1$. If   there exists at least one $k_i$ not equal to $3$ and $G$ is  an $F_{\ell}$-free  graph  of  of  sufficiently large  order $n$, then
$$e(G)\leq  \binom{p}{2}+ p(n-p)+c,$$
where $c=1$ if all $k_i$ are odd and $c=0$ otherwise. Moreover,
if $c=1$ then  equality holds if and only if $G=S^+_{n,p}$. Otherwise, the equality holds if and only if  $G=S_{n,p}$.
\end{Lemma}
Bushaw and Kettle \cite{Bushaw2011}, and Yuan and Zhang \cite{Yuan2017} determined the Tur\'{a}n number for linear forests with all path of order equal to $3$.  Denote by $M_n=\lfloor\frac{n}{2}\rfloor K_2\cup (n-2\lfloor\frac{n}{2}\rfloor K_2)$.
\begin{Lemma}\label{lem2.3}\cite{Bushaw2011,Yuan2017}
Let  $G$ be a  $\ell P_3$-free graph of order $n$. Then 
\[e(G)\leq\left\{
\begin{array}{llll}
 \vspace{1mm}
 \binom{n}{2},&& \mbox{for $n<3\ell$};\\
   \vspace{1mm}
  \binom{3\ell-1}{2}+\big\lfloor\frac{n-3\ell+1}{2}\big\rfloor,&& \mbox{for $3\ell\leq n<5\ell-1$};\\
       \vspace{1mm}
 \binom{3\ell-1}{2}+\ell, && \mbox{for $n=5\ell-1$};\\
      \vspace{1mm}
 \binom{\ell-1}{2}+(n-\ell+1)(\ell-1)+\big\lfloor\frac{n-\ell+1}{2}\big\rfloor, && \mbox{for $n>5\ell-1$.}
\end{array}\right.
\]
Moreover, (i) If $n <3\ell$, then equality holds if and only if $G = K_n$;\\
(ii) If $3\ell\leq n < 5\ell-1$, then  equality holds if and only if  $G= K_{3\ell-1}\cup M_{n-3\ell+1}$;\\
(iii) If $n = 5\ell - 1$, then equality holds if and only if  $G= K_{3\ell-1}\cup M_{2\ell}$ or $G= F_{5\ell-1,\ell-1}$;\\
(iv) If $n > 5\ell - 1$, then equality holds if and only if $G = F_{n,\ell-1}$.\end{Lemma}

Recently, Chen \emph{et al.} \cite{Chen2022} determined the bipartite Tur\'{a}n number for linear forests. For the purpose of this paper, we only use the  upper bound of the bipartite Tur\'{a}n number for linear forests, which can be deduced from \cite[Theorem~1.5]{Chen2022} directly. For more results on linear forests, readers are referred to  \cite{Zhu2021,Zhu2022}.

\begin{Lemma}\label{lem2.4}
 Let $F_{\ell}=\bigcup_{i=1}^{\ell} P_{k_i}$ be a linear forest with   ${\ell}\geq 1$ and
$k_1\geq \cdots \geq k_{\ell}\geq2$ and $ p=\sum _{i=1}^{\ell} \lfloor k_i/2 \rfloor-1$.
If  $n$ is sufficiently larger with  comparing to $p$ and $m$,  then $\ex(m,n;F_{\ell})<pn$.
\end{Lemma}
%
%




\begin{Lemma}\label{Lem2.5}\cite{Chen2022-2}
Let  $\alpha \in[0,1)$, $p\geq1$, and  $n\geq p$. Then
\begin{eqnarray*}
 \rho_\alpha(S_{n,p}) &=& \frac{\alpha n+p-1+\sqrt{(\alpha n+p-1)^2-4p(2\alpha-1)n-4p(p-(p+1)\alpha)}}{2}\\
   &\geq&  \alpha(n-1)+(1-\alpha)(p-1).
\end{eqnarray*}
  In particular, if $n\geq \frac{(2p-1)^2}{2\alpha^2}-\frac{8p^2-2p-1}{2\alpha}+2p(p+1)$ and $\alpha \in(0,1)$, then
$\rho_\alpha(S_{n,p})\geq \alpha n+\frac{2p-1-(2p+1)\alpha}{2\alpha}$.
\end{Lemma}

\begin{Lemma}\label{Lem2.6}
Let  $\alpha \in[0,1)$, $p\geq1$, and  $n\geq p+3$. Then
$\rho_\alpha(S^+_{n,p})$ is equal to the largest root of $f(x)=0$,
where $f(x)=x^3 - (p + n\alpha + p\alpha)x^2 + (p - np + n\alpha - 2 p\alpha + p^2 + np\alpha^2 + 2np\alpha - 1)x + p^2\alpha^2 - 2p + p^3\alpha^2 + np + 5 p\alpha - 2 \alpha^2 + p^2\alpha - p^3\alpha - p^2 + np^2\alpha - 2np^2\alpha^2 - 2np\alpha$.
\end{Lemma}

\begin{Proof}
Set for short $\rho_\alpha =\rho_\alpha (S^+_{n,p})$ and let ${\bf x}$  be a positive
eigenvector corresponding to $\rho_\alpha$. By symmetry,  all vertices with degree $n-1$, $p+1$, and $p$ have the same eigenvector entries, denoted by $x_1$, $x_2$ and $x_3$, respectively.
By eigenequations of $A_\alpha(G)$, we have
\begin{eqnarray*}
   (\rho_\alpha-\alpha(n-1)-(1-\alpha)(p-1))x_1&=& (1-\alpha)(2x_2+(n-p-2)x_3) \\
 (\rho_\alpha-\alpha(p+1))x_2 &=& (1-\alpha)(px_1+x_2)\\
  (\rho_\alpha-\alpha p)x_3 &=& (1-\alpha)px_1.
\end{eqnarray*}
Then it is easy to check that $\rho_\alpha(G)$ is the largest root of $f(x)=0$.
\end{Proof}
\begin{Lemma}\label{Lem2.7}
Let  $\alpha \in[0,1)$, $p\geq1$, and  $n\geq p+3$. If $n-p$ is even, then
 $\rho_\alpha(F_{n,p})=\frac{\alpha n+p+\sqrt{(\alpha n-3p)^2+4(p-\alpha)n+4(p+2)(\alpha p-3p+5)}}{2}$.
 If $n-p$ is odd, then  $\rho_\alpha(F_{n,p})$ is the largest root of $f(x)-p(\alpha-1)^2(n-p-3)=0$, where $f(x)$ is the cubic polynomial in Lemma~\ref{Lem2.6}.
\end{Lemma}

\begin{Proof}
Set for short $\rho_\alpha =\rho_\alpha (F_{n,p})$ and let ${\bf x}$  be a positive
eigenvector corresponding to $\rho_\alpha$.
If $n-p$ is even, then  by symmetry,  all vertices with degree $n-1$ and $p+1$ have the same eigenvector entries, denoted by $x_1$ and $x_2$, respectively.
 By eigenequations of $A_\alpha(G)$, we have
\begin{eqnarray*}
   (\rho_\alpha-\alpha(n-1)-(1-\alpha)(p-1))x_1&=& (1-\alpha)(n-p)x_2 \\
 (\rho_\alpha-\alpha(p+1))x_2 &=& (1-\alpha)(px_1+x_2).
\end{eqnarray*}
Then it is easy to check that $\rho_\alpha$ is the largest root of $ h(x)=0$, where $h(x)=x^2 - ( p + n\alpha)x - ((p - 1)(\alpha - 1) - \alpha(n - 1))(\alpha p + 1) - p(n - p)(\alpha - 1)^2$, i.e.,
$\rho_\alpha=\frac{\alpha n+p+\sqrt{(\alpha n-3p)^2+4(p-\alpha)n+4(p+2)(\alpha p-3p+5)}}{2}$.
If $n-p$ is odd, then  by symmetry, all vertices with degree $n-1$, $p+1$, and $p$ have the same eigenvector entries, denoted by $x_1$, $x_2$ and $x_3$, respectively.
By eigenequations of $A_\alpha(G)$, we have
\begin{eqnarray*}
   (\rho_\alpha-\alpha(n-1)-(1-\alpha)(p-1))x_1&=& (1-\alpha)((n-p-1)x_2+x_3) \\
 (\rho_\alpha-\alpha(p+1))x_2 &=& (1-\alpha)(px_1+x_2)\\
  (\rho_\alpha-\alpha p)x_3 &=& (1-\alpha)px_1.
\end{eqnarray*}
Then it is easy to check that $\rho_\alpha$ is the largest root of $f(x)-p(\alpha-1)^2(n-p-3)=0$.
\end{Proof}

\section{Proof of Theorem~\ref{linear forest}}

We first  prove the following structural lemma for $F_{\ell}$-free connected  graphs of sufficiently large order $n$.
\begin{Lemma}\label{lem3.1}
Let  $F_{\ell}=\cup_{i=1}^{\ell} P_{k_i}$ be a linear forest with   $\ell\geq1$ and
$k_1\geq \cdots \geq k_{\ell}\geq2$ and $ p=\sum _{i=1}^{\ell} \lfloor\frac{k_i}{2}\rfloor-1\geq 1$.  If  $G$ is  an $F_{\ell}$-free connected  graph  of sufficiently large order $n$ and $\rho_{\alpha}(G)\geq \rho_{\alpha}(S_{n,p})$ with $\alpha \in(0,1)$, then there  exists a sufficiently small $\delta$ depending on $n$ and a set $A$ of size $p$  such that $| N[A]|>(1-\delta)n$.
\end{Lemma}

\begin{Proof}
Let  $\rho_{\alpha}=\rho_{\alpha}(G)$ and $\mathbf x$ be a positive eigenvector corresponding to $\rho_{\alpha}$ such that $\mathbf x$ has maximum entry equal to $1$. Choose $w\in V(G)$ with $x_w=1$.
 Let $r=\sum_{i=1}^{\ell} k_{i}$. Set $L=\{v\in V(G): x_v> \epsilon\} $ and $S=\{v\in V(G): x_v\leq \epsilon\} $, where
 $$\epsilon=\max\bigg\{\sqrt{\frac{2r(2p+1)}{\alpha(n-1)}},\sqrt{\frac{2(1-\alpha)(p+1)}{\alpha^2(n-1)}}\bigg\}<\frac{1}{2(3p+2)}.$$
By Lemma~\ref{Lem2.5},
\begin{eqnarray}\label{I2}
  \rho_\alpha &\geq& \rho_\alpha(S_{n,p})\geq \max\bigg\{\alpha n+\frac{2p-1-(2p+1)\alpha}{2\alpha},\alpha (n-1)\bigg\}.
\end{eqnarray}
 By Lemmas~\ref{lem2.1}-\ref{lem2.3},
\begin{eqnarray}\label{I3}
  2e(S)\leq 2e(G)\leq(2p+1)n.
\end{eqnarray}
Next we  prove the following claims.

\vspace{2mm}
{\bf Claim~1.}  $ 2e(L)<\epsilon n$ and $e(L,S)< pn$.
\vspace{2mm}

By eigenequation of $A_\alpha$ on any vertex $u\in L$, we have
$$(\rho_\alpha-\alpha d(u))\epsilon<(\rho_\alpha-\alpha d(u)) x_u=(1-\alpha)\sum_{v\in N(u)} x_v\leq (1-\alpha)d(u),$$
which implies that $$d(u)> \frac{\rho_\alpha\epsilon}{1-\alpha+\alpha\epsilon}>\rho_\alpha\epsilon$$
Thus $$2e(G)=\sum_{u\in V(G)}d(u)\geq \sum_{u\in L }d(u)> |L|\rho_\alpha\epsilon.$$
Combining with (\ref{I2}) and (\ref{I3}), we have
\begin{eqnarray}\label{I0}
|L|<\frac{2e(G)}{\rho_\alpha\epsilon}<\frac{(2p+1)n}{\alpha(n-1)\epsilon}\leq \frac{\epsilon n}{r},
\end{eqnarray}
where the last inequality holds as $\epsilon\geq \sqrt{\frac{2r(2p+1)}{\alpha(n-1)}}$.
Note that the subgraph induced by $L$ is $P_r$-free. By Lemma~\ref{lem2.1} and (\ref{I0}),  $$2e(L)\leq( r-2)|L|<  \frac{(r-2)\epsilon n}{r}< \epsilon n.$$
In addition, since $|S|$ is sufficiently large with comparing to  $|L|$ which follows from (\ref{I0}), we have $e(L,S)< pn$ by Lemma~\ref{lem2.4}. This proves Claim~1.

\vspace{2mm}
{\bf Claim~2.}
\vspace{2mm}
Let $u\in L$. Then for any $u\in L$, we have
$ d(u)> (1-\epsilon)n$.

By eigenequation of $A_\alpha$, we have
\begin{eqnarray*}
  \rho_\alpha \sum_{v\in V(G)}x_v& =& \sum_{v\in V(G)} \rho_\alpha x_v=\sum_{v\in V(G)}\bigg(\alpha d(v)x_v+(1-\alpha)\sum_{z\in N(v)}x_z\bigg)\\
  &=&\alpha\sum_{v\in V(G)} d(v)x_v+(1-\alpha)\sum_{vz\in E(G)}(x_v+x_z)\\
 &=&\sum\limits_{vz\in E(G)}(x_v+x_z)\\
 &=&\sum\limits_{vz\in E(L)}(x_v+x_z)+\sum\limits_{vz\in E(S)}(x_v+x_z)+\sum\limits_{vz\in E(L,S)}(x_v+x_z)\\
  &\leq& 2e(L)+2\epsilon e(S)+(1+\epsilon)e(L,S).
\end{eqnarray*}
Combining with Claim~1, (\ref{I2}) and (\ref{I3}), we have
\begin{equation}\label{I1}
\begin{aligned}
 \sum_{v\in V(G)}x_v&\leq \frac{ 2e(L)+2\epsilon e(S)+(1+\epsilon)e(L,S)}{\rho_\alpha}\\
 &\leq \frac{ \epsilon n+\epsilon (2p+1)n+(1+\epsilon)pn}{\alpha (n-1)}\\
 &<\frac{ (p+1)n}{\alpha (n-1)}.
 \end{aligned}
\end{equation}
In addition, by eigenequations of $A_\alpha$ on  $u$, we have
\begin{eqnarray}\label{I4}
(\rho_\alpha-\alpha d(u)) x_u&=&(1-\alpha)\sum_{v\in N(u)} x_v\leq (1-\alpha)\sum_{v\in V(G)}x_v.
\end{eqnarray}
By  (\ref{I2}), (\ref{I1}), and (\ref{I4}), we have
\begin{eqnarray*}
d(u) &\geq&\frac{\rho_\alpha}{\alpha}-\frac{(1-\alpha)\sum_{v\in V(G)}x_v}{\alpha x_u}
  > \frac{\alpha(n-1)}{\alpha}- \frac{(1-\alpha) (p+1)n}{\alpha^2 (n-1) x_u}\\
  & > &n-1- \frac{(1-\alpha) (p+1)n}{\alpha^2 (n-1) \epsilon}
  \geq n-\frac{\epsilon n}{2}-\frac{\epsilon n}{2} \\
  &=&(1-\epsilon)n,
\end{eqnarray*}
where the last second inequality holds as $\epsilon\geq \sqrt{\frac{2(1-\alpha)(p+1)}{\alpha^2(n-1)}}$.
This proves Claim~2.

\vspace{2mm}

{\bf Claim~3.} $|L|=p$ and $|N[L]|> (1-p\epsilon )n$.

\vspace{2mm}
If $|L|\geq p+1$, then let $B\subseteq L$ be a vertex set of size $p+1$ whose element is arbitrarily chosen from $L$.  By Claim~2, $|N[B]|>(1-(p+1)\epsilon )n\geq 2p+2$. Then
$G$ contains  a subgraph isomorphic to $K_{p+1,2p+2}$, which implies that $G$ contains a copy of $F_{\ell}$, which is a contradiction.
If $|L|\leq p-1$, then  

$$e(L,S)\leq (p-1)(n-p+1).$$ Further, by the definition of $L$ and $S$, we have
\begin{eqnarray}\label{I5}
\sum\limits_{uv\in E(L,S)}(x_u+x_v)\leq(1+\epsilon)e(L,S)\leq (1+\epsilon)(p-1)(n-p+1).
\end{eqnarray}
By eigenequations of $A_\alpha$ on  $w$, we have
$$\rho_\alpha-\alpha d(w)=(\rho_\alpha-\alpha d(w)) x_w=(1-\alpha)\sum_{v\in N(w)} x_v.$$
Multiplying both sides of the above inequality by $\rho_\alpha$, we have
\begin{eqnarray*}
  &&\rho_\alpha(\rho_\alpha-\alpha d(w))\\
  &=& (1-\alpha)\sum_{v\in N(w)} \rho_\alpha x_v \\
   &=& (1-\alpha)\sum_{v\in N(w)} \bigg(\alpha d(v)x_v+(1-\alpha)\sum_{u\in N(v)} x_u\bigg)\\
 &=& (1-\alpha)\sum_{vw\in E(G)} \alpha d(v)x_v+(1-\alpha)^2\sum_{v\in N(w)}\sum_{u\in N(v)} x_u\\
 &\leq&(1-\alpha)\bigg(\sum_{v\in V(G)} \alpha d(v)x_v-\alpha d(w)\bigg)+(1-\alpha)^2\sum_{uv\in E(G)}( x_u+x_v)-\\
 &&(1-\alpha)^2\sum_{v\in N(w)}x_v\\
&=&\alpha(1-\alpha)\sum_{uv\in E(G)}( x_u+x_v)-\alpha(1-\alpha)d(w)+(1-\alpha)^2\sum_{uv\in E(G)}( x_u+x_v)-\\
&&(1-\alpha)(\rho_\alpha-\alpha d(w))\\
&=&(1-\alpha)\sum_{uv\in E(G)}( x_u+x_v)-(1-\alpha)\rho_\alpha,
\end{eqnarray*}
which implies that $$\sum_{uv\in E(G)}( x_u+x_v)\geq \frac{\rho_\alpha(\rho_\alpha+1-\alpha-\alpha d(w))}{1-\alpha}.$$
On the other hand,
\begin{eqnarray*}
  \sum\limits_{uv\in E(G)}(x_u+x_v)&=&\sum\limits_{uv\in E(L,S)}(x_u+x_v)+\sum\limits_{uv\in E(S)}(x_u+x_v) +\sum\limits_{uv\in E(L)}(x_u+x_v)\\
   &\leq& \sum\limits_{uv\in E(L,S)}(x_u+x_v)+ 2\epsilon e(S)+2e(L).
\end{eqnarray*}
 Combining with (\ref{I2})-(\ref{I3}), Claim~1 and $d(w)\leq n-1$, we have
\begin{eqnarray*}
 &&\sum\limits_{uv\in E(L,S)}(x_u+x_v) \\
 &\geq&  \frac{\rho_\alpha(\rho_\alpha+1-\alpha-\alpha d(w))}{1-\alpha}- 2\epsilon e(S)-2e(L)\\
       &\geq& \bigg(\frac{\alpha n}{1-\alpha}+\frac{2p-1-(2p+1)\alpha}{2\alpha(1-\alpha)}\bigg)\bigg(\alpha n+\frac{2p-1-(2p+1)\alpha}{2\alpha}+1-\alpha-\alpha n+\alpha\bigg)-\\
       &&   \epsilon (2p+1)n-\epsilon n\\
       &=&  \bigg(\frac{\alpha n}{1-\alpha}+\frac{2p-1-(2p+1)\alpha}{2\alpha(1-\alpha)}\bigg)\frac{(2p-1)(1-\alpha)}{2\alpha}-\epsilon (2p+2)n\\
        &=&\bigg(p-\frac{1}{2}-\epsilon(2p+2)\bigg)n+\frac{(2p-1)^2-(2p-1)(2p+1)\alpha}{4\alpha^2}\\
       &>&\bigg(p-\frac{1}{2}-(2p+3)\epsilon\bigg) n   \end{eqnarray*}
Combining with (\ref{I5}), we have
\begin{eqnarray*}
(1+\epsilon)(p-1)(n-p+1)
> \bigg(p-\frac{1}{2}-(2p+3)\epsilon\bigg) n,
\end{eqnarray*}
which implies that $$\bigg(\frac{1}{2}-(3p+2)\epsilon\bigg)n< -(p-1)^2(1+\epsilon), $$ which is a contradiction as $(3p+2)\epsilon< \frac{1}{2}$.
By the discussion above, we have  $|L|=p$. By Claim~2, $|N[L]|> (1-p\epsilon )n$. Let $A=L$ and $\delta=p\epsilon$. Then the result follows.
\end{Proof}

\vspace{3mm}

\noindent{\bf Proof of Theorem~\ref{linear forest}.}  Let $G$ be an $F_{\ell}$-free  graph of   order  $n$ with the maximum $A_\alpha$-spectral radius. Let $\rho_\alpha=\rho_\alpha(G)$ and $\mathbf x$ be a positive eigenvector corresponding to $\rho_{\alpha}$. Let
\[H_{n,p}=\left\{
\begin{array}{llll}
 \vspace{1mm}
 S_{n,p},&& \mbox{there  exists  at least one even $k_i$};\\
   \vspace{1mm}
 S^+_{n,p},&& \mbox{all $k_i$  are  odd and there exists one $k_i>3$};\\
       \vspace{1mm}
F_{n,p}, && \mbox{$k_1=\cdots=k_{\ell}= 3$.}
\end{array}\right.
\] By Lemmas~\ref{Lem2.5}-\ref{Lem2.7}, it suffices to prove that $G=H_{n,p}$. Next we consider the following two cases.

{\bf Case~1.}  $G$ is connected. Since $S_{n,p}$ is $F_{\ell}$-free, $\rho_\alpha\geq \rho_\alpha(S_{n,p})$ by the extremality of $G$. By Lemma~\ref{lem3.1}, there  exists a sufficiently small $\delta$ depending on $n$ and a set $A$ of size $p$  such that $|N[A]|>(1-\delta)n$.
Next let $B= N[A]$ and $R=V(G)\backslash (A\cup B)$.

(i) First suppose that there  exists  at least one even $k_i$.  Then $R\cup B$ induces isolated vertices, otherwise $G$ contains a copy of $F_{\ell}$ as $G$ is connected.
If $R$ is not empty, then $G$ is a proper subgraph of $S_{n,p}$, which implies that $\rho_\alpha<\rho_\alpha(S_{n,p})$, a contradiction.   Hence $R$  is empty. By the extremality of $G$ again, $G=S_{n,p}$.

(ii) Next suppose that all $k_i$  are  odd and there exists one $k_i>3$.  Since $S^+_{n,p}$ is $F_{\ell}$-free, we have $\rho_\alpha\geq \rho_\alpha(S^+_{n,p})$ by the extremality of $G$. Further, $G[B]$ contains at most one edge, otherwise  the subgraph induced by $A\cup B$  contains a copy of $F_{\ell}$, a contradiction.

We first assume that $G[B]$ contains precisely one edge, denoted by  $e$. Then $G[R\cup B]$ also contains precisely one edge $e$, otherwise $G$ contains a copy of $F_{\ell}$  as $G$ is connected. If $R$ is not empty, then $G$ is a proper subgraph of $S^+_{n,p}$, which implies that $\rho_\alpha<\rho_\alpha(S^+_{n,p})$, a contradiction.  Thus $R$  is empty.  By the extremality of $G$ again, $G=S^+_{n,p}$.

Next assume that $G[B]$ consists of  isolated vertices. Then $G[R\cup B]$ is $P_3$-free, otherwise $G$ contains a copy of $F_{\ell}$   as $G$ is connected.  Hence $G[R\cup B]$ consists of
independent edges and isolated vertices. i.e., $v_1v_2,\ldots, v_{2t-1}v_{2t}, w_1,\ldots, w_r$.
If  $t\leq1$ then
  $G$ is a proper subgraph of $ S^+_{n,p}$ and thus
 $\rho_\alpha<\rho_\alpha(S^+_{n,p})$, a contradiction.
So  $t\geq2$. If $p=1$, then $F_{\ell}=P_5$. In this situation, since $G$ is connected, we have $r=0$ and every edge $v_{2i-1}v_{2i}$ has one end in $B$ and the other in $R$ for all $1\leq i\leq t$. Then  $G$ contains a copy of $P_5$, a contradiction. So $p\geq2$. Since $G$ is $F_{\ell}$-free,  there is no edge with one end in $B$  and the other end in $R$.  This implies that $v_{2i-1}v_{2i}\in E(G[R])$ for all $1\leq i\leq t$.
 Likewise,  for all $1\leq j\leq 2t$, we see that any $v_j$   is adjacent to precisely one vertex in $A$,  denoted by $u$.
 We next show  that $\rho_\alpha<\rho_\alpha(S^+_{n,p})$.
 Choose a vertex $v\in A\backslash \{u\}$. 
 From the eigenequations  of $\rho_\alpha(G)$, we have  for $1\leq i\leq t$,
 $$(\rho_\alpha-2\alpha)x_{v_{2i-1}}=(1-\alpha)(x_{v_{2i}}+x_u), \quad(\rho_\alpha-2\alpha)x_{v_{2i}}=(1-\alpha)(x_{v_{2i-1}}+x_u),$$
 $$(\rho_\alpha- \alpha d(v))x_v=(1-\alpha)\sum \limits_{w\in N(v)}x_w>(1-\alpha)x_u.$$
Note that $d(v)> (1-\delta)n>1+\frac{1}{\alpha}$.  Then $$x_{v_{2i-1}}=x_{v_{2i}}=\frac{(1-\alpha)x_u}{\rho_\alpha-1-\alpha}, \qquad x_v>\frac{(1-\alpha)x_u}{\rho_\alpha- \alpha d(v)}>x_{v_{2i-1}}.$$
Let $G'$ be a graph obtained from $G$ by deleting all edges in $\{v_{2i-1}v_{2i}: 1\leq i\leq t\}$ and
adding all edges in  $\{vv_i: 1\leq i\leq 2t\}$.  Obviously, $G'$ is a proper subgraph of $ S_{n,p}$ and
$$\rho_\alpha(G')< \rho_\alpha(S_{n,p})<\rho_\alpha(S^+_{n,p}).$$ On the other hand,
\begin{eqnarray*}
 &&\rho_\alpha(G')-\rho_\alpha(G) \\
 &\geq &\frac{\sum \limits_{1\leq i \leq 2t}(\alpha x_v^2+2(1-\alpha)x_vx_{v_{i}}+\alpha x_{v_i}^2)-\sum\limits_{1\leq i \leq t} (\alpha x_{v_{2i-1}}^2+2(1-\alpha)x_{v_{2i-1}}x_{v_{2i}}+\alpha x_{v_{2i}}^2)}{{\bf x}^\mathrm{T}{\bf x}}\\
 &>&0,
\end{eqnarray*}
implying that $\rho_\alpha<\rho_\alpha(G')$. Thus $\rho_\alpha<\rho_\alpha(S^+_{n,p})$, a contradiction.

(iii) Suppose that $k_1=\cdots=k_{\ell}= 3$.  Since $G$ is $F_{\ell}$-free, we have $G[R\cup B]$ is $P_3$-free. Hence $G[R\cup B]$ consists of independent edges and isolated vertices. If $R$ is not  empty, then $G$ is a proper subgraph of $F_{n,p}$,  which implies that $\rho_\alpha<\rho_\alpha(F_{n,p})$, a contradiction. Thus $R$ is empty.  By the extremality of $G$, $G=F_{n,p}$.

\vspace{2mm}
{\bf Case~2.}  $G$ is not connected.  Since $S_{n,p}$ is $F_{\ell}$-free,
$$\rho_\alpha\geq \rho_\alpha(S_{n,p})\geq \alpha n+\frac{2p-1-(2p+1)\alpha}{2\alpha}.$$ Let $G_1$ be a component of $G$  such that $\rho_\alpha(G_1)=\rho_\alpha$. Set $n_1=|V(G_1)|$.
Then
\begin{eqnarray*}
  n_1-1 \geq  \rho_\alpha(G_1)=\rho_\alpha\geq \alpha n+\frac{2p-1-(2p+1)\alpha}{2\alpha},
\end{eqnarray*}
which implies that
$n_1$ is also sufficiently large.
 By Case~1, $\rho_\alpha(G_1)\leq \rho_\alpha(H_{n_1,p})$.
 Since  $\rho_\alpha(G_1)=\rho_\alpha$ and  $\rho_\alpha(H_{n_1,p})< \rho_\alpha(H_{n,p})$,  we have
 $$\rho_\alpha< \rho_\alpha(H_{n,p}).$$ On the other hand, since $H_{n,p}$ is $F_{\ell}$-free, we have $\rho_\alpha\geq \rho_\alpha(H_{n,p})$, a contradiction. This completes the proof.
\QEDB


\begin{thebibliography}{12}

\bibitem{Bushaw2011} N. Bushaw, N. Kettle, Tur\'{a}n numbers of multiple paths and equibipartite forests, Combin. Probab. Comput. 20 (2011) 837--853.

\bibitem{Chen2019} Y.Y. Chen, D.Li, Z.W. Wang, J.X. Meng, $A_\alpha$-spectral radius of the second power of a graph, Appl. Math. Comput. 359 (2019) 418--425.
      \bibitem{Chen2020}  M.-Z. Chen, A-M. Liu, X.-D. Zhang, The signless Laplacian spectral radius of graphs with forbidding linear forests, Linear Algebra Appl. 591 (2020) 25--43.

  \bibitem{Chen2022}  M.-Z. Chen, Ning Wang, Long-Tu Yuan, X.-D. Zhang, The bipartite Tur\'{a}n number and spectral extremum for linear forests, arXiv:2201.00453.

   \bibitem{Chen2022-2}   M.-Z. Chen, A-M. Liu, X.-D. Zhang, Spectral extremal results on the $\alpha$-index of graphs without minors and star forests, arXiv:2204.00181.

 \bibitem{Erdos1959}  P. Erd\H{o}s,  T. Gallai, On maximal paths and circuits of graphs. Acta Math. Acad. Sci. Hungar. 10 (1959) 337--356.


\bibitem{Guo2022} S.G. Guo,   R. Zhang, The sharp upper bounds on the $A_\alpha$-spectral radius of $C_4$-free graphs and Halin graphs, Graphs Combin. 38 (2022)  Paper No. 19, 13 pp.
\bibitem{He2021}C.X. He, W.Y. Wang, Y.Y. Li, L.L Liu, Some Nordhaus-Gaddum type results of $A_\alpha$-eigenvalues of weighted graphs, Appl. Math. Comput. 393 (2021) 125761.



\bibitem{Lidicky2013} B. Lidicky, H. Liu, C. Palmer, On the Tur\'{a}n number of forests, Electron. J. Combin. 20 (2) (2013) 62.

 \bibitem{Lin2019} H.Q.  Lin,  X.G. Liu, J. Xue,  Graphs determined by their $A_\alpha$-spectra, Discrete Math. 342 (2019) 441--450.
%

\bibitem{Nikiforov2017} V. Nikiforov,  Merging the $A$- and $Q$-spectral theories. Appl. Anal. Discrete Math. 11 (1) (2017) 81--107.


\bibitem{Tian2021}  G.X. Tian, Y.X. Chen, S.Y. Cui, The extremal $\alpha$-index of graphs with no 4-cycle and 5-cycle, Linear Algebra Appl. 619 (2021) 160--175.
%
 \bibitem{Yu2019} Z.Q. Yu, L.Y. Kang, L.L. Liu, E.F. Shan, The extremal $\alpha$-index of outerplanar and planar graphs, Appl. Math. Comput., 343 (2019) 90--99.
\bibitem{Yuan2017}  L.-T. Yuan, X.-D. Zhang, The Tur\'{a}n number of disjoint copies of paths, Discrete Math. 340 (2) (2017) 132--139.
%

\bibitem{Zhu2021} X.T. Zhu, F.F. Zhang, Y.J. Chen, Generalized Tur\'{a}n number of even linear forests, Graphs Combin. (2021) 1437--1449.
 \bibitem{Zhu2022} X.T. Zhu,  Y.J. Chen, Generalized Tur\'{a}n number for linear forests, Discrete Math. 345 (2022) 112997.
\end{thebibliography}
\end{document}